\documentclass[10pt]{amsart}
\input{epsf.tex}

\usepackage{latexsym,amsfonts,amssymb,epsfig,verbatim}
\usepackage{amsmath, latexsym, graphics, textcomp,psfrag}
\usepackage{float}
\usepackage{amsthm}

\usepackage[all,2cell,dvips]{xy}

\newcommand{\nc}{\newcommand}
\nc{\dmo}{\DeclareMathOperator}
\nc{\nt}{\newtheorem}

\nc{\ds}{\displaystyle}
\nc{\ens}{\ensuremath}

\theoremstyle{plain}
\nt{thm}{Theorem}[section]
\nt*{main}{Main Theorem}
\nt{lem}[thm]{Lemma}
\nt{cor}[thm]{Corollary}
\nt*{cor*}{Corollary}
\nt{prop}[thm]{Proposition}
\nt*{q}{Question}
\nt{fact}{Fact}
\nt*{claim}{Claim}
\nt{conjecture}[thm]{Conjecture}

\theoremstyle{remark}

\newtheorem{definition-remark}[thm]{Definition-Remark}

\newtheorem*{definition*}{Definition}
\newtheorem*{remark*}{Remark}

\def\S{{\Sigma}}

\def\C{{\mathcal C}}

\nc{\N}{\mathcal{N}}

\def\S{{\mathcal S}}

\dmo{\out}{Out}
\dmo{\aut}{Aut}
\dmo{\gl}{GL}
\dmo{\SL}{SL}
\dmo{\sy}{Sp}
\dmo{\mcg}{Mod}

\dmo{\stu}{\tilde{St}}
\dmo{\st}{St}
\dmo{\lku}{\tilde{Lk}}
\dmo{\lk}{Lk}
\dmo{\dlk}{Lk_<}
\dmo{\cdlk}{Lk_\ll}

\nc{\T}{\mathcal{T}}
\nc{\tm}{\mathcal{M}}
\nc{\X}{\mathcal{X}}
\nc{\Y}{\mathcal{Y}}
\nc{\K}{\mathcal{K}}
\nc{\IA}{\textrm{IA}}

\nc{\sg}{\Sigma_g}
\nc{\modg}{\mcg(\sg)}
\nc{\mods}{\mcg(S)}
\nc{\spz}{\sy_{2g}(\Z)}

\nc{\Z}{\mathbb{Z}}
\nc{\R}{\mathbb{R}}
\nc{\Q}{\mathbb{Q}}
\nc{\I}{\mathcal{I}}
\nc{\W}{W}

\nc{\outf}{\out(F_n)}
\nc{\glz}{\gl_n(\Z)}
\nc{\outfn}{\out(F_n)}
\nc{\glnz}{\gl_n(\Z)}
\nc{\glzw}{W \backslash \! \glz}

\nc{\kpo}{K(\T_n,1)}

\nc{\rt}{\tilde{\rho}}
\nc{\strt}{\st(\rt)}
\nc{\strto}{\st(\rt_1)}
\nc{\strtt}{\st(\rt_2)}
\nc{\strtop}{\st(\rt_1')}
\nc{\strttp}{\st(\rt_2')}
\nc{\lkrt}{\lk(\rt)}
\nc{\str}{\st(\rho)}
\nc{\lkr}{\lk(\rho)}
\nc{\dlkr}{\dlk(\rho)}
\nc{\cdlkr}{\cdlk(\rho)}

\nc{\li}{L_i}
\nc{\lit}{\tilde\li}

\nc{\seg}{\mathcal{S}}
\nc{\lseg}{\underline{\mathcal{S}}}

\nc{\slice}{\Sigma_I}
\nc{\bzt}{\ens{{\mathbb Z}^2}}

\nc{\norm}[1]{\|#1\|}
\nc{\bignorm}[1]{\left\|#1\right\|}

\nc{\tbt}[4]{\ensuremath{\left(\begin{array}{rr} #1 & #2 \\ #3 & #4\end{array}\right)}}

\nc{\twbo}[2]{\ensuremath{\left(\begin{array}{c} #1 \\ #2 \\  \end{array}\right)}}

\nc{\tbo}[3]{\ensuremath{\left(\begin{array}{c} #1 \\ #2 \\ #3 \end{array}\right)}}

\nc{\ttt}[9]{\ensuremath{\left(\begin{array}{rrr} #1 & #2 & #3 \\ #4 & #5 & #6 \\ #7 & #8 & #9  \end{array}\right)}}

\nc{\p}[1]{\paragraph{{\bf #1}}}

\nc{\pics}[3]{\epsfysize=#3 cm \begin{figure}[htb]
\center{\leavevmode \epsfbox{#1.eps}} \caption{#2}  \label{#1pic}
\end{figure}}

\nc{\beqn}{\begin{equation}}
\nc{\eeqn}{\end{equation}}

\nc{\bl}{ \begin{list}{$\cdot$}{
\setlength{\leftmargin}{.5in}
\setlength{\rightmargin}{.5in}
\setlength{\parsep}{0.5ex plus .2ex minus 0ex}
\setlength{\itemsep}{0.2ex plus 0.2ex minus 0ex}
}
}
\nc{\blwide}{ \begin{list}{$\cdot$}{
\setlength{\leftmargin}{.25in}
\setlength{\rightmargin}{.25in}
\setlength{\parsep}{0.5ex plus .2ex minus 0ex}
\setlength{\itemsep}{0.2ex plus 0.2ex minus 0ex}
}
}

\nc{\el}{\end{list}}

\nc{\bpf}{\begin{proof}}
\nc{\epf}{\end{proof}}

\nc{\si}{\sigma}
\nc{\ep}{\epsilon}
\nc{\al}{\alpha}
\nc{\be}{\beta}
\nc{\Ga}{\Gamma}
\nc{\la}{\lambda}

\nc{\margin}[1]{\marginpar{\scriptsize #1}}

\nc{\set}[1]{\{#1\}}

\nc{\bc}{}

\begin{document}

\title[Dimension of the Torelli group for $\out(F_n)$]{Dimension of the Torelli group for \boldmath$\out(F_n)$}

\author{Mladen Bestvina}
\author{Kai-Uwe Bux}
\author{Dan Margalit}

\address{Mladen Bestvina: Department of Mathematics\\ University of Utah\\ 155 S 1440 East \\ Salt Lake City, UT 84112-0090}
\email{bestvina@math.utah.edu}

\address{Kai-Uwe Bux: Department of Mathematics\\ University of Virginia\\ Kerchof
  Hall 229\\ Charlottesville, VA 22903-4137}
\email{kb2ue@virginia.edu}

\address{Dan Margalit: Department of Mathematics\\ University of Utah\\ 155 S 1440 East \\ Salt Lake City, UT 84112-0090}
\email{margalit@math.utah.edu}

\thanks{The first and third authors gratefully acknowledge support by
  the National Science Foundation.}



\maketitle

\begin{center}\today\end{center}

\begin{abstract}
Let $\T_n$ be the kernel of the natural map $\outf \to \glz$.  We
use combinatorial Morse theory to prove that $\T_n$ has an Eilenberg--MacLane space which is $(2n-4)$-dimensional and that
$H_{2n-4}(\T_n,\Z)$ is not finitely generated ($n \geq 3$).  In particular, this
recovers the result of Krsti\'c--McCool that $\T_3$ is not finitely
presented.  We also give a new proof of the fact, due to Magnus, that
$\T_n$ is finitely generated.
\end{abstract}


\section{Introduction}

There is a natural homomorphism from $\out(F_n)$, the group of outer
automorphisms of the free group on $n$ generators, to $\gl_n(\Z)$,
given by abelianizing the free group $F_n$.  It is a theorem of
Nielsen that this map is surjective \cite{jn}.  We call its kernel the
\emph{Torelli subgroup} of $\out(F_n)$, and we denote it by $\T_n$:

\[ 1 \to \T_n \to \out(F_n) \to \gl_n(\Z) \to 1 \]

\begin{main} \label{m1}
For $n \geq 3$, we have:
\begin{enumerate}
\item $\T_n$ has a $(2n-4)$-dimensional Eilenberg--MacLane space.
\item $H_{2n-4}(\T_n,\Z)$ is infinitely generated.
\item \label{fg} $\T_n$ is finitely generated.
\end{enumerate}
\end{main}

Part~(\ref{fg}) of the main theorem is due to Magnus; we give our own
proof in Section~\ref{section:finite generation}.  We remark that
$\T_1$ is obviously trivial and $\T_2$ is trivial by a classical
result of Nielsen \cite{jn} (we give a new proof of the latter fact in
Section~\ref{section:finite generation}).

The group $\T_n$, like any torsion free subgroup of $\out(F_n)$, acts
freely on the spine for outer space (see Section~\ref{section:kpo}),
and therefore has an Eilenberg--MacLane space of dimension $2n-3$, the dimension of this
spine.  Our theorem improves this upper bound on the dimension and
shows that $2n-4$ is sharp.

When $n=3$, we obtain that $H_2(\T_3,\Z)$ is not finitely generated,
and this immediately implies the result of Krsti\'c--McCool that
$\T_3$ is not finitely presented \cite{km}.

\p{Historical background} The question of whether $H_k(\T_n,\Z)$ is finitely generated, for
various values of $k$ and $n$, is a long standing problem with few
solutions.  This question was explicitly asked by Vogtmann in her
survey article \cite{kv}.  We now give a brief history of related
results, all of which are recovered by our main theorem.

Nielsen proved in 1924 that $\T_3$ is finitely generated \cite{jn}.  Ten years later, Magnus proved
that $\T_n$ is finitely generated for every $n$ \cite{wm}.

Smillie--Vogtmann proved in 1987 that, if $2 < n < 100$ or $n > 2$ is
even, then $H_\star(\T_n,\Z)$ is not finitely
generated \cite{sv1} \cite{sv2}.  Their method is to consider the
rational Euler characteristics of the groups in the short exact
sequence defining $\T_n$ (see \cite{kv}).

The Krsti\'c--McCool result that $\T_3$
is not finitely presented was proven in 1997, via completely algebraic
methods \cite{km}.
It is a general fact that if the second homology of a group is not
finitely generated, then the group is not finitely presented.

\p{Large abelian subgroups} It follows from the second part of the
main theorem that the first part is sharp; i.e., $\T_n$ does not have
an Eilenberg--MacLane space of dimension less than $2n-4$.  A simpler proof that the
cohomological dimension of $\T_n$ is at least $2n-4$ is to simply
exhibit an embedding of $\Z^{2n-4}$ into $\T_n$.  There is a subgroup
$\Z^{2n-4} \cong G < \T_n$ consisting of elements with representative automorphisms
given by:
\[
\begin{array}{rcl}
x_1 & \mapsto & x_1 \\
x_2 & \mapsto & x_2 \\
x_3 & \mapsto & [x_1,x_2]^{p_3}x_3[x_1,x_2]^{q_3} \\
& \vdots & \\
x_n & \mapsto & [x_1,x_2]^{p_n}x_n[x_1,x_2]^{q_n}
\end{array}
\]
for varying $p_i$ and $q_i$ (the $x_i$ are generators for $F_n$).

In Section~\ref{section:infinite generation}, we prove that specific
conjugates of $G$ represent independent classes in
$H_{2n-4}(\T_n,\Z)$, thus proving the second part of the main theorem.
These conjugates are exactly the generators of $H_{2n-4}(\tm_n,\Z)$,
where $\tm_n$, called the ``toy model'', is a particularly simple
subcomplex of the Eilenberg--MacLane space $\Y_n$ defined in
Section~\ref{section:kpo}.  In Section~\ref{section:infinite
generation}, we prove that $H_{2n-4}(\tm_n,\Z)$ injects into
$H_{2n-4}(\Y_n,\Z)$.  Since the homology of $\tm_n$ is not finitely
generated in any dimension greater than 1, we are led to the following
question:

\begin{q}
\label{question:middle homology}
Does $H_\ast(\tm_n,\Z)$ inject into $H_\ast(\Y_n,\Z)$?
\end{q}

\p{Mapping class groups.} The term ``Torelli group'' comes from the
theory of mapping class groups.  Let $\sg$ be a closed surface of
genus $g \geq 1$.  The \emph{mapping class group} of $\sg$, denoted
$\modg$, is the group of isotopy classes of orientation preserving
homeomorphisms of $\sg$.  The \emph{Torelli group}, $\I_g$, is the
subgroup of $\modg$ acting trivially on the homology of $\sg$.  As
$\modg$ acts on $H_1(\sg,\Z)$ by symplectic automorphisms, $\I_g$ is
defined by:
\[ 1 \to \I_g \to \modg \to \spz \to 1 \]
It is a classical theorem of Dehn, Nielsen, and Baer that the natural
map $\modg \to \out(\pi_1(\sg))$ is an isomorphism.  In this sense
$\T_n$ is the direct analog of $\I_g$.

Our (lack of) knowledge of the finiteness properties of $\I_g$ mirrors
that for $\T_n$.  Using the fact that $\mcg(\Sigma_1) \cong \SL_2(\Z) =
\sy_2(\Z)$, it is
obvious that $\I_1$ is trivial.  Johnson showed in 1983 that $\I_g$ is
finitely generated for $g \geq 3$ \cite{dj1}.  In 1986,
McCullough--Miller showed that $\I_2$ is not finitely generated
\cite{mcm}, and Mess improved on this in 1992 by showing that $\I_2$
is a free group of infinite rank \cite{gm}.  At the same time, Mess
further showed that $H_3(\I_3,\Z)$ is not finitely generated.  In
Kirby's problem list, Mess asked about finiteness properties in higher
genus \cite{rk}.

\p{Automorphisms vs. outer automorphisms.}  Strictly speaking, Magnus
and Krsti\'c--McCool study the group $\K_n$, by which we mean the
kernel of $\aut(F_n) \to \gl_n(\Z)$, where $\aut(F_n)$ is the
automorphism group of the free group.  By considering the short
exact sequence
\[ 1 \to F_n \to \K_n \to \T_n \to 1 \]
we see that $\K_n$ is finitely generated if and only if $\T_n$ is
finitely generated.  Moreover, it follows from our main theorem and the spectral
sequence associated to this short exact sequence
that $H_{2n-3}(\K_n,\Z)$ is not finitely generated and if $k$ is the
smallest index so that $H_k(\T_n,\Z)$ is not finitely generated, then
$H_k(\K_n,\Z)$ is not finitely generated.

From a topological point of view, $\T_n$ is the more natural group to
study.

In the literature, $\T_n$ is sometimes denoted by
$\IA_n$ for ``identity on abelianization'' (see, e.g. \cite{kv}).
However, since Krsti\'c--McCool use $\IA_n$ to denote the kernel of $\aut(F_n)
\to \gl_n(\Z)$, we avoid this notation to eliminate the confusion.
The notation $\K_n$ comes from Magnus \cite{wm}.

\p{Acknowledgements.}
We would like to thank Bob Bell, Mikhail Gromov, Jon McCammond, and Kevin
Wortman for helpful conversations.  We are especially grateful to Karen Vogtmann for
explaining her unpublished work.


\section{An Eilenberg--MacLane space}
\label{section:kpo}

In Section~\ref{section:spine}, we recall
the definition of Culler--Vogtmann's spine for Outer space.  Then, in
Section~\ref{section:quotient}, we describe the quotient of this space by
$\T_n$.  This quotient is a $(2n-3)$-dimensional Eilenberg--MacLane space for $\T_n$.

A \emph{rose} is a graph with one vertex.  The \emph{standard rose} in
rank $n$, denoted $R_n$, is a particular rose which is fixed once and
for all.  We denote the standard generators of $F_n \cong \pi_1(R_n)$
by $x_1, \dots, x_n$.


\subsection{Spine for Outer space} \label{section:spine}
Culler--Vogtmann introduced the \emph{spine for Outer space}, which we
denote by $\X_n$, as a tool for studying $\out(F_n)$ \cite{cv}.  This
is a simplicial complex defined in terms of marked graphs.

A \emph{marked graph} is a pair $(\Ga,g)$, where $\Ga$ is a finite
metric graph (1-dimensional cell complex with a metric) with no
separating edges and no vertices of valence less than 3 and $g:R_n \to
\Ga$ is a homotopy equivalence ($g$ is called the \emph{marking}).  We
say that two marked graphs $(\Ga,g)$ and $(\Ga',g')$ are
\emph{equivalent} if $g' \circ g^{-1}$ is homotopic to an isometry,
where $g^{-1}$ is any homotopy inverse of $g$.  We will denote the
equivalence class $[(\Ga,g)]$ by $(\Ga,g)$.

The vertices of $\X_n$ are equivalence classes of marked graphs where
all edges have length 1.  A set of vertices
\[ \set{(\Ga_1,g_1), \dots, (\Ga_k,g_k)} \]
is said to span a simplex if $\Ga_{i+1}$ is obtained from $\Ga_i$ by
collapsing a forest in $\Ga_i$, and $g_{i+1}$ is the marking obtained
from $g_i$ via this operation.

We can think of arbitrary points of $\X_n$ as marked metric graphs:
for instance, as we move along an edge between two vertices in $\X_n$,
the length of some edge in the corresponding graphs (more generally,
the lengths of the edges in a forest) varies between 0 and 1.

There is a natural right action of $\out(F_n)$ on $\X_n$.  Namely,
given $\phi \in \out(F_n)$ and $(\Ga,g) \in \X_n$, the action is
given by:
\[ (\Ga,g) \cdot \phi = (\Ga,g \circ \phi) \]
(here we are using the fact that every element $\phi$ of $\out(F_n)$
can be realized by a homotopy equivalence $R_n \to R_n$, also denoted
$\phi$, uniquely up to homotopy).

Culler--Vogtmann proved the following result \cite{cv}:

\begin{thm}
\label{cv thm}
For $n \geq 2$, the space $\X_n$ is contractible.
\end{thm}

This theorem has the consequence that the virtual cohomological
dimension of $\out(F_n)$ is equal to $2n-3$, the dimension of
$\X_n$.

The \emph{star of a rose} in $\X_n$ is the union of the closed
simplices containing the vertex corresponding to a rose.  The key idea
for Theorem~\ref{cv thm} is to think of $\X_n$ as the union of stars of
vertices corresponding to marked roses.  We take an analogous approach
in this paper.


\subsection{The quotient}\label{section:quotient}

Baumslag--Taylor proved that $\T_n$ is torsion free \cite{bt}.  We
also know that the action of $\T_n$ on $\X_n$ is free: by the
definition of the action, point stabilizers correspond to graph
isometries, and isometries act nontrivially on homology.
Finally, the action is simplicial, and so it follows
that the quotient of $\X_n$ by $\T_n$ is an Eilenberg--MacLane space for $\T_n$:
\[ \Y_n = \X_n/\T_n \]

\p{Homology markings.} Since $\out(F_n)$ identifies every pair of
isometric graphs of $\X_n$, points of $\Y_n$ can be thought of as
equivalence classes of pairs $(\Ga,g)$, where $\Ga$ is a metric graph
(as before), and $g$ is a \emph{homology marking}; that is, $g$ is an
equivalence class of homotopy equivalences $R_n \to \Ga$, where two
homotopy equivalences are equivalent if (up to isometries of $\Ga$)
they induce the same map $H_1(R_n,\Z) \to H_1(\Ga,\Z)$.

Via the marking $g$, we can think of the (oriented) edges of $\Ga$ as
elements of $H^1(R_n) \cong \Z^n$ (if $e$ is an edge and $x$ is a
simplicial 1-chain, then $e(x)$ is the number of times $e$
appears in $x$).  As such, if we think of the generators $x_1,
\dots, x_n$ of $\pi_1(R_n)$ as elements of $H_1(R_n,\Z)$, then we can
\emph{label} each oriented edge $a$ of $\Ga$ by the corresponding row
vector:
\[ (a(g(x_1)), \dots, a(g(x_n))) \]
where $a(g(x_i))$ is the number of times
$g(x_i)$ runs over $a$ homologically, with sign.

In this way, a point of $\Y_n$ is given by a labelled graph, and two
such graphs represent the same point in $\Y_n$ if and only if there is
a label preserving graph isomorphism between them (i.e. if the map
induces the identity on cohomology).  See Figure~\ref{figure:not
  frontier} for an example of a labelled graph.  We remark that this
example exhibits the fact that $\Y_n$ is not a simplicial
complex---there are two edge collapses (and hence two edges in $\Y_n$)
taking this point to the rose with the identity marking.

\begin{figure}[htb]
\psfrag{x1}{$(1,0,0)$}
\psfrag{x2}{$(0,1,0)$}
\psfrag{x3}{$(0,0,1)$}
\centerline{\includegraphics[scale=.35]{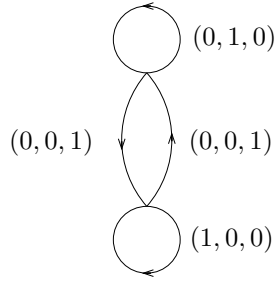}}
\caption{An example of a labelled graph.}
\label{figure:not frontier}
\end{figure}

When convenient, we will confuse the points of $\Y_n$ with the
corresponding marked graphs.

We will make use of the following generalities about marked graphs in $\Y_n$:

\begin{prop}
\label{lemma1}
Let $(\Ga,g)$ be a marked graph.
\begin{enumerate}
\item \label{labels dont change} If an edge of $\Ga$ is collapsed, the labels of the
  remaining edges do not change.
\item \label{parallel}Any two edges of $\Ga$ with the same label (up
  to sign) are \emph{parallel}
  in the sense that the union of their interiors disconnects $\Ga$.
\item \label{switch condition} The sum of the labels of the (oriented) edges coming into a
vertex of $\Ga$ is equal to the sum of the labels of the edges leaving the
vertex.
\end{enumerate}
\end{prop}

We leave the proofs to the reader.

\p{Roses.}  Let $\Ga$ be a rose with edges $a_1, \dots, a_n$, and let
$g:R_n \to \Ga$ be a homology marking.  Up to isometries of $\Ga$, the marking $g$ gives an element of
$\gl_n(\Z)$, called the \emph{marking matrix}; the rows are exactly the labels of the edges.

Since all edges have length 1, the isometry group of $\Ga$ is generated
by swapping edges and by reversing the orientations of edges; the
former operation has the effect of switching rows of the matrix, and
the latter corresponds to changing signs of rows.  Thus, in this case,
$(\Ga,g)$ gives rise to an element of $\W\backslash\gl_n(\Z)$, where
$\W=\W_n$ is the signed permutation subgroup of $\gl_n(\Z)$, acting on
the left.  In fact, this gives a bijection between roses in $\Y_n$ and
elements of $\W\backslash\gl_n(\Z)$, as $\out(F_n)$ acts transitively
on the roses of $\X_n$.

\medskip

The right action
of $\out(F_n)$ on $\X_n$ descends to a right action of $\gl_n(\Z)$ on
$\Y_n$.  In particular, the action on roses is given by the
right action of $\gl_n(\Z)$ on $\W\backslash\gl_n(\Z)$.


\section{Stars of roses}
\label{section:star of a rose}
\label{section:frontier}

As with $\X_n$, we like to
think of the quotient $\Y_n$ as the union of stars of roses. By
definition, the \emph{star of a rose} in $\Y_n$ is the image of the
star of a rose in $\X_n$.  Thus, it consists of graphs which can be
collapsed to a particular rose.
We now discuss some of the basic properties of the star of a rose.


\subsection{Labels in the star of a rose} We will need several observations
about the behavior of labels in the star of a rose.  The proofs of the
various parts of the propositions are straightforward and are left to
the reader.  In each of the statements, let $\rho$ be a rose in $\Y_n$
represented by a marked graph $(\Ga,g)$.  Say that its edges $a_1,
\dots, a_n$ are labelled by $v_1, \dots, v_n \in \Z^n$.

\begin{prop}
\label{lemma2}
If $(\Ga',g')$ is a marked graph in $\str$, we have:
\begin{enumerate}
\item \label{vi} For each $i$, there is an edge of $\Ga'$ labelled $\pm v_i$.
\item \label{forest}The edges not labelled $\pm v_i$ form a forest.
\item \label{circles}If each edge of $\Ga'$ is labelled $\pm v_i$,
  then the union of edges labelled $\pm v_i$ (for any particular $i$)
  is a topological circle (see, e.g., Figure~\ref{figure:not frontier}).
\item \label{coefficients} The label of any edge of $\Ga'$ is of the form
\[ \sum_{i=1}^n k_i v_i \]
where $k_i \in \{-1,0,1\}$.
\end{enumerate}
\end{prop}

We have the following converse to the first two parts of the previous proposition:

\begin{prop}
\label{converse}
If $(\Ga',g')$ is a marked graph which has, for each $i$, at least one
edge of length 1 labelled $\pm v_i$, then $(\Ga',g')$ is in
$\str$.
\end{prop}

We also have a criterion for when a marked graph is in the frontier of
the star of a rose:

\begin{prop}
\label{lemma:frontier condition}
A marked graph $(\Ga',g')$ in $\str$ is in the frontier
of $\str$ if and only if it has at least one edge of length 1 whose
label is not $\pm v_i$ for any $i$.  In this case, the given label is
a label for some rose whose star contains $(\Ga',g')$.
\end{prop}


\subsection{Ideal edges}
\label{section:ideal edges}

Let $\rho=(\Ga,g)$ be a rose whose edges $a_1, \dots, a_n$ are labelled $v_1, \dots, v_n$,
as above.  An \emph{ideal edge} is any formal sum:
\[ \sum_{i=1}^n k_i a_i \]
where $k_i \in \{-1,0,1\}$, and at least two of the $k_i$ are nonzero.
An ideal edge is a ``direction'' in $\str$ in the following sense:
for any ideal edge, we can find a marked graph $(\Ga',g')$ in the
frontier of $\str$ where one of the edges of $(\Ga',g')$ has the label
$\sum k_i v_i$.  If a marked graph in $\str$ has an edge of length 1
with label $\sum k_i v_i$, we say that the marked graph
\emph{realizes} the ideal edge $\sum k_i a_i$.

\begin{lem}
\label{lem:1 edge}
Given any ideal edge for a particular rose, there is a 1-edge blowup
of $\rho$ in the frontier of $\str$ which realizes that ideal edge.
\end{lem}

\begin{figure}[H]
\psfrag{v1}{$v_1$}
\psfrag{v2}{$v_2$}
\psfrag{v3}{$v_3$}
\psfrag{v4}{$v_4$}
\psfrag{v5}{$v_5$}
\psfrag{v1v3v4}{$v_1-v_3+v_4$}
\centerline{\includegraphics[scale=.5]{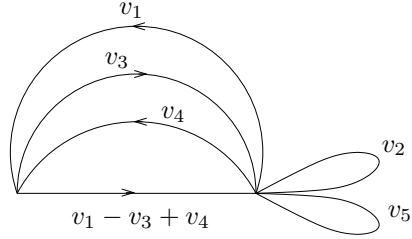}}
\caption{A 1-edge blowup realizing the ideal edge $a_1-a_3+a_4$.}
\label{figure:ideal edge}
\end{figure}

The lemma is proven by example.  See Figure~\ref{figure:ideal edge}
for a picture of a 1-edge blowup realizing the ideal edge
$a_1-a_3+a_4$ in rank 5 (apply Proposition~\ref{lemma1}(\ref{switch
condition})).  Also, we see that there are many graphs satisfying the
conclusion of the lemma---the other graphs are obtained by moving the
loops labelled $v_2$ and $v_5$ arbitrarily around the graph.

Our notion of an ideal edge is simply the homological version of the
ideal edges of Culler--Vogtmann \cite{cv}.

We say that an ideal edge $\iota'$ is \emph{subordinate} to the ideal
edge $\iota = \sum k_i a_i$ if $\iota'$ is obtained by
changing some of the $k_i$ to zero.
A \emph{2-letter ideal edge} is an ideal edge of
the form $k_{i}a_{i} + k_{j}a_{j}$.  Two ideal edges are said
to be \emph{opposite} if one can be obtained from the other by
changing the sign of exactly one coefficient.  The following facts are
used in Section~\ref{section:finite generation}:

\begin{lem}
\label{lemma:subordinate}
Let $\rho=(\Ga,g)$ be a rose whose edges $a_1, \dots, a_n$ are
labelled $v_1, \dots, v_n$.  Suppose that $\iota$ and $\iota'$ are
ideal edges and that either:
\begin{enumerate}
\item $\iota'$ is subordinate to $\iota$, or
\item $\iota$ and $\iota'$ are 2-letter ideal edges which are not opposite.
\end{enumerate}
In either case, there is a marked graph $(\Ga',g')$ in $\str$ which
simultaneously realizes $\iota$ and $\iota'$.
\end{lem}

\bpf

\begin{figure}[H]
\psfrag{v1}{$v_1$}
\psfrag{v2}{$v_2$}
\psfrag{v3}{$v_3$}
\psfrag{v4}{$v_4$}
\psfrag{v1v2}{$v_1+v_2$}
\psfrag{v2v3}{$v_2+v_3$}
\psfrag{v3v4}{$v_3+v_4$}
\psfrag{v1v2v3v4}{$v_1+v_2+v_3+v_4$}
\psfrag{-v1v2}{$-v_1-v_2$}
\centerline{\includegraphics[scale=.5]{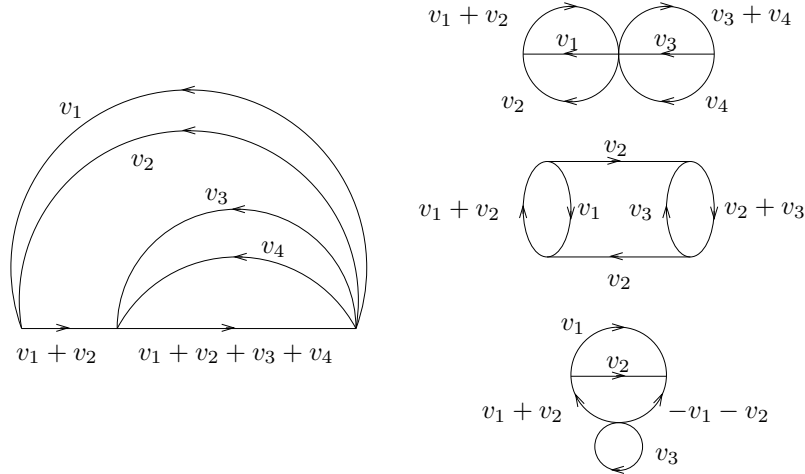}}
\caption{Marked graphs simultaneously realizing subordinate ideal
  edges (left) and 2-letter ideal edges which are not opposite
  (right).}
\label{figure:simultaneous}
\end{figure}

In each case, we can explicitly describe the desired graph.  If
$\iota'$ is subordinate to the ideal edge $\iota = \sum k_i
a_i$, we start with a 1-edge blowup realizing $\iota$
(Lemma~\ref{lem:1 edge}), and then blow
up another edge to separate the edges which appear in $\iota'$ from
those which do not.  Figure~\ref{figure:simultaneous} (left hand side)
demonstrates this for $\iota=v_1+v_2+v_3+v_4$ and $\iota' = v_1+v_2$
in rank 4.

For the case of two 2-letter ideal edges which are not opposite,
without loss of generality it suffices to demonstrate marked graphs
which simultaneously realize $v_1+v_2$ with $-v_1-v_2$, $v_2+v_3$, or
$v_3+v_4$ (the arbitrary case is obtained by renaming/reorienting
edges and by attaching extra 1-cells to any vertex).  See
Figure~\ref{figure:simultaneous} (right hand side) for a
demonstration.  One can use Proposition~\ref{lemma1}(\ref{switch
  condition}) to verify the labels.
\epf

The reader may verify that opposite ideal edges are never
simultaneously realized.

We remark that, in the framework established by Culler--Vogtmann, one
can think of this lemma in terms of compatibility of partitions, in
which case the proof is immediate; see
\cite{cv}.


\subsection{Homotopy type} In the remainder of this section, we
prove that the star of any rose retracts onto the subcomplex consisting
of ``cactus graphs'', and this subcomplex is homeomorphic to a union
of $(n-2)$-tori.

We define a \emph{rank $n$ cactus graph} inductively as follows.  A
rank 1 cactus graph is a graph with 1 vertex and 1 edge (i.e. a circle
with a distinguished point).  In general, a rank $n$ cactus graph is
obtained by gluing a rank 1 cactus graph to a rank $n-1$ cactus graph
along the vertex of the rank 1 cactus graph.  The set of vertices of
the new graph is the union of the sets of vertices of the original two
graphs.
We note that a rank $n$ cactus graph has exactly $n$ embedded
circles, and every edge belongs to exactly one embedded circle
(Figure~\ref{figure:not frontier} is an example).

Let $\C(\rho)$ denote the space of cactus graphs in $\str$.
Given any $\rho'$, there is a canonical homeomorphism $\C(\rho) \to
\C(\rho')$, once we choose orderings of the edges of $\rho$ and
$\rho'$.  Thus, we can unambiguously use $\C_n$ to denote the space of
cactus graphs in the star of a rose in rank $n$.

In the remainder, assume that $\rho=(\Ga,g)$ is a rose in $\Y_n$ with
edges $a_1, \dots, a_n$ labelled $v_1, \dots, v_n$.

\begin{lem}
\label{lemma:retraction to cactus graphs}
$\str$ strongly deformation retracts
onto $\C(\rho)$.
\end{lem}

\bpf

For every marked graph in $\str$, the set of edges
whose label is not $\pm v_i$ is a forest
(Proposition~\ref{lemma2}(\ref{forest})).  We perform a strong deformation
retraction of $\str$ by shrinking the edges of each such forest in
each marked graph in $\str$.

Consider any marked graph $(\Ga',g')$ in the image of the retraction.
By Propositions~\ref{lemma2}(\ref{vi}) and~\ref{lemma2}(\ref{circles}),
there is a circle of edges labelled $\pm v_i$ for each $i$.  We
consider the ``dual graph'' obtained by assigning a vertex to each
such circle (the \emph{circle vertices}) and each intersection point
(the \emph{point vertices}) and we connect a point vertex to a circle
vertex if the point is contained in the circle.  It follows from
Proposition~\ref{lemma1}(\ref{parallel}) that this graph is a tree,
and hence $(\Ga',g')$ is a cactus graph.
\epf

\begin{cor}
\label{star dim}
For $n \geq 2$, the star of any rose $\str$ in $\Y_n$ is homotopy equivalent to a complex of dimension $n-2$.
\end{cor}

\bpf

By the definition of cactus graphs, we can see that the dimension
increases with slope 1 with respect to dimension, starting at $n=2$.
Since $\C_2$ is a point, $\C_n$ is a complex of dimension
$n-2$.  An application of Lemma~\ref{lemma:retraction to cactus
graphs} completes the proof.
\epf

We can filter $\C(\rho)$ by subsets according to the
number of vertices in the cactus graphs:
\[ \set{\rho} = V_0 \subset V_1 \subset \cdots \subset V_{n-2} = \C(\rho) \]
Each $V_i$ consists of cactus graphs with $i-1$ vertices.

Our goal now is to give a generating set for $\pi_1(\C(\rho))$.
Since $V_2$ is simple to understand, the following proposition will
make it easy to do this.

\begin{prop}
\label{prop:cactus space}
There is a cell structure on $\C(\rho)$ so that the $i$-skeleton is
exactly $V_i$.
\end{prop}

\bpf

We proceed inductively.  The 0-skeleton is one point $V_0 = \set{\rho}$.

Let $i > 0$.  Any marked graph $(\Ga',g')$ in $V_i - V_{i-1}$ lies in
a unique $i$-cell $C$.  For each $i$, let $k_i$ be the number of
edges of $\Ga'$ labelled $\pm v_i$.  If we reparameterize so that the sum of the
lengths of the edges of $\Ga'$ labelled $\pm v_i$ is 1, then we get a
$(k_i-1)$-simplex for each $i$, and $C$ is the product of these simplices.

The boundary of $C$ is the set of points where some edge is
assigned length 0.  Clearly, $\partial C \subset V_{i-1}$. and so
the proposition follows.
\epf

For the remainder of this section, we use the cell structure given by
Proposition~\ref{prop:cactus space}, which is different from the cell
structure inherited from $\Y_n$.

Let $V_1^1$ be the subset of $V_1$ consisting of graphs with a vertex
of valence 4 and a vertex of valence $2n-2$ (i.e. only a single loop
is ``travelling'' around another loop).  We will see in
Section~\ref{section:finite generation} that the obvious generators
for $\pi_1(V_1^1)$ correspond to one of the two types of Magnus
generators for $\T_n$.

\begin{prop}
\label{cor:gens for cactus space}
The subcomplex $V_1^1$ contains a generating set for $\pi_1(\C(\rho))$.
\end{prop}

\bpf

First, each of $\C_1$ and $\C_2$ is a single point.  In
rank 3, $V_1^1=V_1$.  Thus, in all of these cases, the proposition is
vacuously true.  For the remainder, assume $n \geq 4$.

As per Proposition~\ref{prop:cactus space}, $V_1$ can be thought of as
the 1-skeleton of the cell complex $\C(\rho)$.  This subcomplex has
1 vertex (the rose $\rho$) and an edge for each combinatorial type of labelled
graph with 2 vertices.  We now need to show that any such \emph{standard loop} $\alpha$ in
$V_1$ can be written in $\pi_1(\C(\rho))$ as a product of loops in
$V_1^1$.  Our strategy is to show that $V_2$ is a union of 2-tori and
that $\pi_1(V_1^1)$ surjects onto $\pi_1(V_2)$.

\begin{figure}[ht]
\includegraphics[scale=.35]{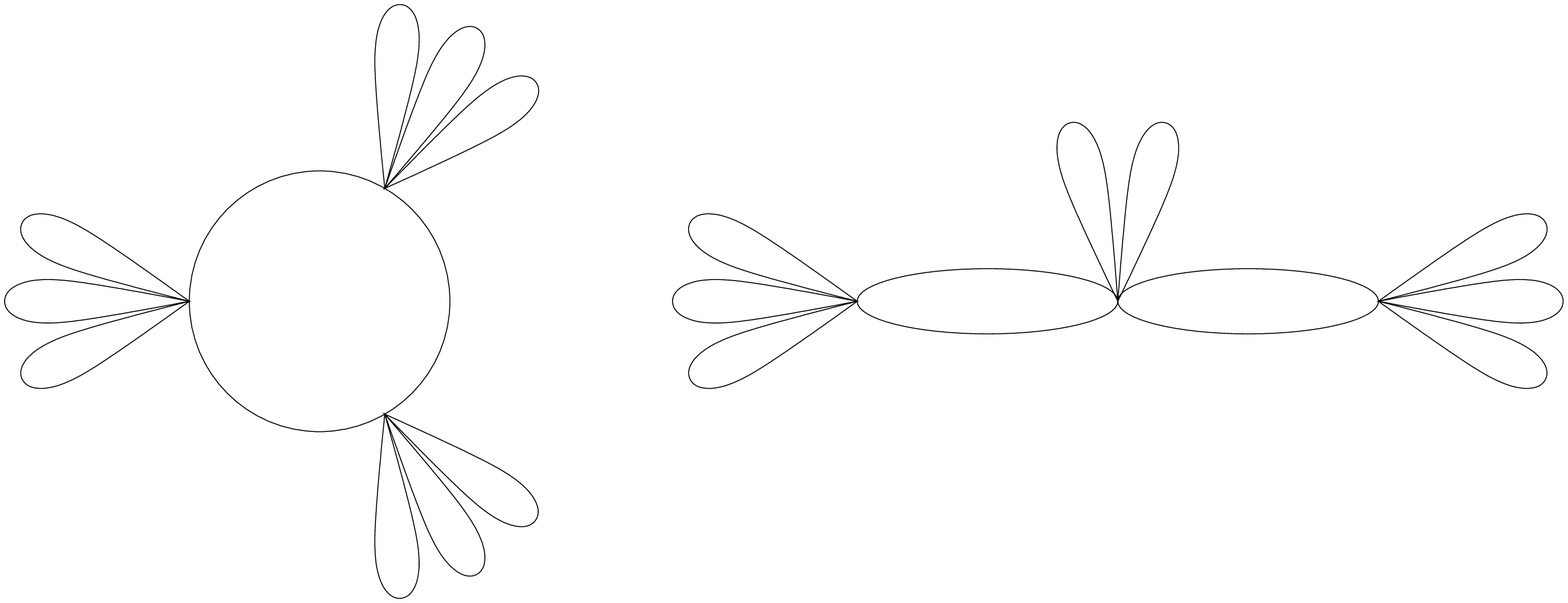}
\caption{Two types of graphs in $V_2$.}\label{figure:v3}
\end{figure}

If $(\Ga',g')$ is a point of $V_2-V_1$, then there are two
possibilities: the three vertices of $\Ga'$ either lie on the same
circle or they do not; see Figure~\ref{figure:v3}.  If they do all lie
on some ``central circle'', then we obtain a 2-torus by fixing one
vertex and letting the other two vertices ``move around'' the central
circle (really we are changing lengths so as to give the appearance of
this motion).  In the other case, there are two central circles.  By
fixing the middle intersection point and letting the other two
intersection points move around the respective circles, we again see a
torus.

Consider a standard loop $\al$ of $V_1$.  At an interior point of $\alpha$,
there is a central circle with two vertices, and the two vertices have
valence, say, $p=p(\al)$ and $q=q(\al)$.  By definition of $V_1$, we
have that $p$ and $q$ are even and at least $4$; say $p \leq q$.  We thus
have a filtration of $V_1$: $\al$ is in $V_1^k$ if $(p-2)/2
\leq k$.  The number $(p-2)/2$ is the number of loops glued to that vertex,
other than the central circle.

Now, suppose that $\al$ is a standard loop of $V_1^k$ for some $k
\geq 2$.  At any interior point of $\al$, we perform a blowup so that
we end up with a graph in $V_2-V_1$ of the first type (left side of
Figure~\ref{figure:v3}).  Moreover, we choose the blowup so that (at
least) one of the vertices has valence 4.  The fundamental group of the
corresponding torus is generated by a standard loop from $V_1^{k-1}$ and a standard loop
from $V_1^1$, and so $\al$ can be written as a product of such loops.
By induction, $\al$ can be written as a product of loops from $V_1^1$.
\epf

\p{Remark} For completeness, we mention that the entire space
$\C_n$ can be thought of as a union of $(n-2)$-tori, and the
intersection between any two of these tori is a lower dimensional
torus which is a product of diagonals of coordinate subtori.  It is
straightforward to prove this, given what we have already done.
However, we will not need this fact.


\section{Cohomological Dimension}
\label{section:morse function}

We now give the argument for the first part of the main theorem, that
$\Y_n$ is homotopically $(2n-4)$-dimensional.  The basic strategy is to
put an ordering on the stars of roses of $\Y_n$ (we think of the
ordering as a Morse function) and then to glue the stars of roses
together in the prescribed order.  This is in the same spirit as the
proof of Culler--Vogtmann that $\X_n$ is contractible.


\subsection{Morse function}  The ordering on roses will come from an
ordering on matrices.  We start with vectors.
By the \emph{norm} of an element $v=(a_1,\cdots,a_n)$ of $\Z^n$, we mean:
\[ |v|=(|a_1|,\cdots,|a_n|)\in \Z_+^n \]
 where the elements of $\Z_+^n$ are ordered lexicographically.
Consider the matrix:
\[ M = \left(\begin{array}{c} v_1 \\ v_2 \\ \vdots \\ v_n
\end{array}\right) \]
The \emph{norm} of $M$ is:
\[ |M| = (|v_n|,\dots,|v_1|) \in (\Z_+^n)^n \]
where $(\Z_+^n)^n$ has the lexicographic ordering on the $n$ factors.
We say that $M$ is a \emph{standard representative} for an element of
$\W\backslash\gl_n(\Z)$ if $|v_n| < \cdots < |v_1|$ (i.e. if it is a
representative with smallest norm).  Note that two rows of a matrix in
$\gl_n(\Z)$ cannot have the same norm, for otherwise these
two rows would be equal after reducing modulo 2, and the resulting
matrix would not be invertible.

We declare the norm of an element of $\W\backslash\gl_n(\Z)$ to be the
norm of a standard representative, and the norm of a
rose in $\Y_n$ to be the norm of the corresponding element of
$\W\backslash\gl_n(\Z)$.

In what follows, the following fact will be important:

\begin{lem}
\label{lemma:neighbors}
If the stars of two roses
intersect, then the roses have different norms.
\end{lem}

\bpf

If $M$ and $M'$ are marking matrices for neighboring roses, then
$M'=NM$, where each entry of $N$ is either $-1$, $0$, or $+1$ (apply Proposition~\ref{lemma2}(\ref{coefficients})).
Then, if $|M|=|M'|$, it follows that $N$ is the identity modulo 2, and
so $N \in \W$.
\epf

The norm on roses turns the set of roses into a well-ordered set.  We
use this fact without mention in the transfinite induction arguments for Theorem~\ref{thm:dim} and Proposition~\ref{proposition:normal generation}.


\subsection{The induction}
\label{subsection:the induction}

We define an \emph{initial segment} of $\Y_n$ to be a union of stars of
a set of roses that is closed under taking smaller roses (i.e. a
sublevel set of the ``Morse function'' given on stars of roses).  Note
that, in general, an initial segment consists of infinitely many
roses.  If we show that each initial segment is $(2n-4)$-dimensional,
it will follow by transfinite induction that $\Y_n$ has the same
property.

To this end, we define the \emph{descending link} of a rose in $\Y_n$
to be the intersection of its star with the union of all stars of
roses of strictly smaller norm (by Lemma~\ref{lemma:neighbors}, we need not worry about roses of equal norm).  The descending link of a rose $\rho$, denoted $\dlkr$, is a subset of the frontier of its star.  We will prove the following in
Section~\ref{section:descending links}:

\begin{prop}
\label{desc link dim}
For $n \geq 3$, descending links are homotopically $(2n-5)$-dimensional.
\end{prop}

Given this, we can prove the first part of the main theorem:

\begin{thm}
\label{thm:dim}
For $n \geq 3$, the complex $\Y_n$ is homotopy equivalent to a complex of dimension at most $2n-4$.
\end{thm}

\bpf

We proceed by transfinite induction on initial segments.  The base step is Corollary~\ref{star dim}.

Whenever we
glue the star of a rose $\str$ to an initial segment $\S$ in order to
make a new initial segment $\hat \S$, we can think of this as a
diagram of spaces:
\[ \S \leftarrow \dlkr \rightarrow \str \]
By the inductive
hypothesis, $\S$ is homotopy equivalent to a $(2n-4)$-dimensional
space $\S'$.  Denote by $\str'$ the $(n-2)$-complex homotopy equivalent to $\str$ given by Proposition~\ref{star dim}.  By
Proposition~\ref{desc link dim}, the descending link $\dlkr$ is homotopy equivalent to a $(2n-5)$-dimensional space
$\dlkr'$.  We choose maps $\dlkr' \to \S'$ and $\dlkr' \to \str'$ so
that the following diagram commutes up to homotopy:

\begin{figure}[H]
\begin{center}
\scalebox{1}{ \xymatrix{\str' \ar@{<->}[r] & \str \\
\dlkr'  \ar@{->}[d] \ar@{->}[u] \ar@{<->}[r] & \dlkr  \ar@{->}[d] \ar@{->}[u] \\
\S' \ar@{<->}[r] & \S}}
\end{center}
\end{figure}
It follows that the colimit of the diagram of spaces in the left column is homotopy equivalent to the colimit of the diagram of spaces
in the right column (see e.g. \cite[Proposition 4G.1]{ah}).  The
former, call it $\hat \S'$, is $(2n-4)$-dimensional (consider the double mapping cylinder), and
the latter is $\hat \S$.  By construction, the homotopy equivalence $\hat S' \to \hat S$ extends the homotopy equivalence $S' \to S$.

By transfinite induction, we thus build a homotopy model $Z_n$ for
$\Y_n$.  By the inductive construction given above, $Z_n$ has a
filtration by subcomplexes $\set{Z_n^\al}$, each equipped with a
homotopy equivalence $h_\al:Z_n^\al \to \Y_n^\al$ for some some initial
segment $\Y_n^\al$.  What is more, the induced map $h:Z_n \to
\Y_n$, when restricted to $Z_n^\al$, is precisely $h_\al$.  It follows
that $h$ is a homotopy equivalence (see, e.g., the discussion
following \cite[Proposition 4G.1]{ah}).  Since $Z_n$ has dimension at
most $2n-4$ (by construction), we are done.
\epf

\p{Remark.} If one wants to avoid transfinite induction, it is
possible to alter the Morse function so that it is the same locally
(i.e. Proposition~\ref{desc link dim} and its proof do not change)
but the image of the Morse function is order isomorphic to the
positive integers.


\section{Finite generation}
\label{section:finite generation}

In this section, we recall the definition of the Magnus generating set
for $\T_n$, and explain how our point of view recovers the result that
these elements do indeed generate $\T_n$
(Theorem~\ref{proposition:finite generation} below).

Throughout the section (and the appendix), we denote an element $\phi$ of $\out(F_n)$
by:
\[ [\Phi(x_1), \dots, \Phi(x_n)] \]
where $x_1,\dots,x_n$ are the generators of $F_n$, and $\Phi$ is a
representative automorphism for $\phi$.


\subsection{Magnus generators}

Magnus proved that $\T_n$ is generated by:
\[ \begin{array}{rcl}
K_{ik} & = & [x_1, \dots, x_kx_ix_k^{-1}, \dots, x_n] \\
K_{ikl} & = & [x_1, \dots, x_i[x_k,x_l], \dots, x_n]
\end{array}\]
for distinct $i$, $k$, and $l$.

We can see the $K_{ik}$ as loops in the star of a rose in $\Y_n$.
Consider the picture in Figure~\ref{figure:not frontier}.  As
mentioned in Section~\ref{section:quotient}, shrinking
either of the parallel edges gives a path leading to the rose with the
identity marking, and so this is a loop in the star of that rose in
$\Y_n$.  By considering what is happening on the level of homotopy (as
opposed to homology), we see that this loop is exactly $K_{23}$ (see
\cite{cv}).  By attaching more loops at one of the vertices, and
renaming the edges, we see that we can obtain any $K_{ik}$ in the star
of the identity rose.  In the stars of other roses, the
analogously defined loops are conjugates of the $K_{ik}$.  What is
more, we have:

\begin{prop}
\label{proposition:generators of star of a rose}
The fundamental group of the star of the rose with the identity
marking is generated by the $K_{ik}$.
\end{prop}

The proposition follows immediately from the fact that the loops in
the above discussion corresponding to the $K_{ik}$ are exactly the
standard generators for $\pi_1(V_1)$ from Proposition~\ref{cor:gens
for cactus space}.


\subsection{Proof of finite generation}

Our proof that the Magnus generators generate $\pi_1(Y_n) \cong \T_n$
rests on the following two topological facts about descending links
which we prove in Section~\ref{section:descending links}:

\begin{prop}
\label{proposition:descending links are nonempty}
Descending links are nonempty, except for that of the rose with the identity marking.
\end{prop}

\begin{prop}
\label{proposition:descending links are connected}
Descending links are connected.
\end{prop}

Combining Propositions~\ref{proposition:generators of star of a rose},
\ref{proposition:descending links are nonempty},
and~\ref{proposition:descending links are connected} with Van
Kampen's theorem and the transitivity of the action of $\out(F_n)$ on
stars of roses, we see that the fundamental group of any initial segment of $\Y_n$ is normally generated by the $K_{ik}$.  By transfinite induction, we have:

\begin{prop}
\label{proposition:normal generation}
$\T_n$ is normally generated by the $K_{ik}$.
\end{prop}

The group generated by the $K_{ik}$
is not normal in $\out(F_n)$, as any element
of this subgroup is of the form:
\[ [g_1x_1g_1^{-1},g_2x_2g_2^{-1},\dots,g_nx_ng_n^{-1}] \]
Thus, to find a generating set for $\T_n$, we need to add more
elements.

We have the following result of Magnus:

\begin{prop}
\label{magnus}
For any $n$, the group generated by
\[ \set{K_{ik},K_{ikl}: i \neq k < l \neq i} \]
is normal in $\out(F_n)$.
\end{prop}

It is now easy to prove the following, which is the third part of our
main theorem:

\begin{thm}
\label{proposition:finite generation}
$\T_n$ is finitely generated.  In particular, it is generated by $\set{K_{ik},K_{ikl}}$.
\end{thm}

Proposition~\ref{magnus} is also one of the steps in
Magnus's proof that the $K_{ik}$ and $K_{ikl}$ generate $\K_n$
\cite{wm}.  For completeness, we give Magnus's proof of
Proposition~\ref{magnus} in the appendix.


\subsection{Proof that \boldmath$\T_2$ is trivial} Since there are two
ways to blow up a rank 2 rose, it follows that the star of a rose in
$\Y_2$ is homeomorphic to an interval and that the frontier is
homeomorphic to $S^0$.  If we glue the stars of roses together
inductively according to our Morse function as in
Section~\ref{section:morse function}, then at each stage we are gluing
a contractible space (the star of the new rose) to a contractible
space (the previous initial segment is contractible by induction)
along a contractible space (Propositions~\ref{proposition:descending links
are nonempty} and~\ref{proposition:descending links are connected} and
the fact that the frontier is $S^0$).  It follows that each initial
segment, and hence all of $\Y_2$, is contractible; hence, $\T_2=1$.

\begin{figure}[H]
\centerline{\includegraphics{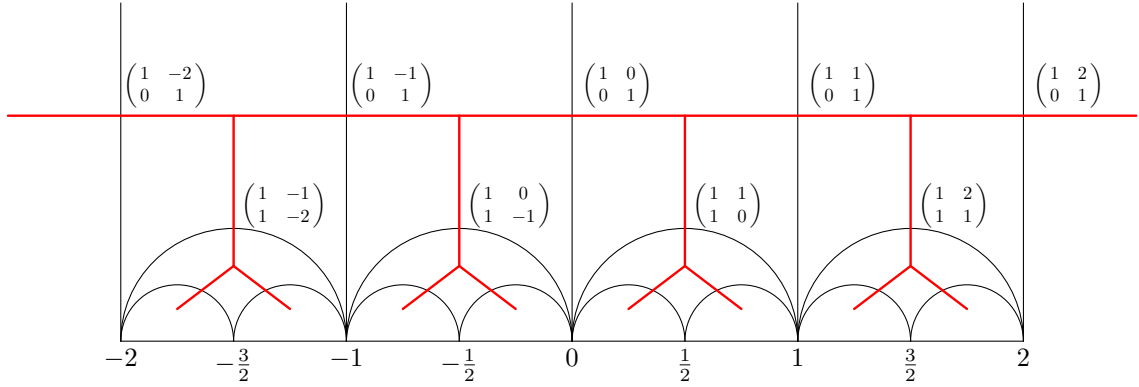}}
\caption{A part of the Farey graph and the dual tree $\Y_2$.}
\label{id18}
\end{figure}

It is more illuminating to draw a diagram of $\X_2=\Y_2$. It is a
tree, with edges representing stars of roses.  This tree is naturally
dual to the classical Farey graph, with the matrix
$\left(\begin{smallmatrix} a&b\\ c&d\end{smallmatrix}\right)$
corresponding to the unordered pair $\{\frac ba,\frac dc\}$. See
Figure \ref{id18}.

\section{Descending links}
\label{section:descending links}

Recall that the descending link $\dlkr$ of a rose $\rho$ is the
intersection of its star with the union of stars of roses of strictly
smaller norm.  The goal of this section is to prove
Propositions~\ref{proposition:descending links are nonempty},
\ref{proposition:descending links are connected}, and~\ref{desc link
dim}, that descending links are nonempty, connected, and homotopically
$(2n-5)$-dimensional.

As in Section~\ref{section:frontier}, let $\rho$ be a rose represented by a marked graph $(\Ga,g)$ whose edges $a_i$ are labelled $v_i$.  We assume the $a_i$ are ordered so that the marking matrix
\[ M = \left(\begin{array}{c} v_1 \\ \vdots \\  v_n \end{array}\right) \]
is a standard representative.


\subsection{Descending ideal edges}

An ideal edge for $\rho$ is called \emph{descending} if any of the corresponding
1-edge blowups (Lemma~\ref{lem:1 edge}) lies in $\dlkr$.  Every edge of a marked graph in $\str$ which is not labelled $\pm v_i$ corresponds to some ideal edge; if the corresponding ideal edge is descending, we may say that the edge is descending.

We now give a criterion for checking whether or not a particular ideal
edge is descending.

\begin{lem}
\label{descending criterion}
Let $\iota = a_1 + a_{i_1} + \cdots + a_{i_m}$ be an ideal edge.  The
following are equivalent:
\begin{enumerate}
\item $\iota$ is descending
\item any of the corresponding 1-edge blowups lies in $\dlkr$
\item all of the corresponding 1-edge blowups lie in $\dlkr$
\item $|v_1 + v_{i_1} + \cdots + v_{i_m}| < |v_1|$
\end{enumerate}

Similarly, $\bar \iota = -a_1 + a_{i_1} + \cdots + a_{i_m}$ is descending if and only if
$|v_1 - (v_{i_1} + \cdots + v_{i_m})| < |v_1|$.
\end{lem}

\bpf

A 1-edge blowup which realizes the ideal edge $\iota$ lies in $m+1$
stars of roses (Proposition~\ref{converse}).  Namely, for each of $v_1, v_{i_1}, \dots v_{i_m}$,
we get a new marking matrix by replacing that vector with:
\[ v_1 + v_{i_1} + \cdots + v_{i_m} \]
and leaving all other row vectors the same.  To see if $\iota$ is
descending, we look at the smallest of these matrices.
We claim that the smallest is:
\[ N = \left(\begin{array}{c} v_1 + v_{i_1} + \cdots + v_{i_m} \\ v_2 \\ \vdots \\  v_n \end{array}\right) \]
Indeed, suppose we had replaced some other row vector, say $v_i$, with
$v_1 + v_{i_1} + \cdots + v_{i_m}$, obtaining a matrix $N'$.  Now,
forgetting the order of the rows, $N$
and $N'$ share $n-1$ rows, and $N$ has the row vector $v_i$ whereas
$N'$ has the row vector $v_1$.  By the assumption that $M$ is a
standard representative, we have $|v_i| < |v_1|$.  Now, if we put
$|N'|$ in standard form, it is easy to find a representative for the
$N$-coset with smaller norm than the standard representative for
$N'$---simply replace the row of $N'$ consisting of $v_1$ with the
vector $v_i$.  The norm of $N$ is less than or equal to the norm of
this representative, so the claim is proven.

Now both directions are easy: if $|v_1 + v_{i_1} + \cdots + v_{i_m}| <
|v_1|$ then $|N|$ is obviously strictly less than $|M|$ (the given
representative has smaller norm) and so $\iota$ is
descending; conversely, if $|v_1 + v_{i_1} + \cdots + v_{i_m}| \geq
|v_1|$, then the given representative is in standard form and
obviously has norm at least $|M|$.  (We remark that the last
inequality must be strict by Lemma~\ref{lemma:neighbors}.)

The second statement follows by symmetry.
\epf

It is not hard to prove a stronger statement than the one given here.  However, the relatively simple result given suffices for our purposes, and the generalities are notationally unpleasant.


\begin{cor}
\label{lemma:description of descending}
A marked graph in $\str$ is in $\dlkr$ if and only if it realizes a descending ideal edge.
\end{cor}

As a consequence of Lemma~\ref{descending criterion}, we see that
there exist pairs of marked graph which can never be simultaneously
descending.

\begin{lem}
\label{forbidden pairs}
If the ideal edge $\iota = a_1 + a_{i_1} + \cdots + a_{i_m}$ is descending
then $\bar \iota = -a_1 +
a_{i_1} + \cdots + a_{i_m}$ is not descending.
\end{lem}

Generalizations of Lemma~\ref{descending criterion} lead to analogous generalizations of the current lemma.

\bpf

To simplify notation, let $w_0=v_1$, $w_1 = v_{i_1}$,
$w_2=v_{i_2}$, etc.  We will denote particular entries in each of
these row vectors by using double indices; i.e., $w_{jk}$ is the
$k^{\mbox{\tiny th}}$ entry of the row vector $w_j$.

Let $k$ be the smallest number so that
\[ |w_{0k} + w_{1k} + w_{2k} + \cdots + w_{mk}| \neq |w_{0k}| \]
Note that there is such a $k$, for otherwise, the
original matrix $M$ would not be invertible (reduce modulo 2).

Applying Lemma~\ref{descending criterion}, we see that $\iota$ is descending if and only if
\begin{equation}\tag{1}
|w_{0k} + w_{1k} + w_{2k} + \cdots + w_{mk}| < |w_{0k}|
\end{equation}
(we are using the minimality of $k$).  It follows that $w_{1k} +
w_{2k} + \cdots + w_{mk}\neq 0$ and that the sign of this sum differs
from that of $w_{0k}$.  Thus, we have:
\begin{equation}\tag{2}
|w_{0k} - (w_{1k} + w_{2k} + \cdots + w_{mk})| > |w_{0k}|
\end{equation}
and so $\bar \iota$ is not descending.  By symmetry, we are done.
\epf


\subsection{Proof of Propositions~\ref{proposition:descending links are nonempty} and~\ref{proposition:descending links are connected}}

As usual, let $\rho$ be a rose represented by
a marked graph $(\Ga,g)$ with edges $a_1, \dots, a_n$ labelled by
$v_1, \dots, v_n$, and assume that the edges are ordered so that the
marking matrix
\[ M = \left(\begin{array}{c} v_1 \\ \vdots \\  v_n \end{array}\right) \]
is a standard representative.

We first give the proof that descending links are nonempty:

\begin{proof}[Proof of Proposition~\ref{proposition:descending
  links are nonempty}]

Let $k$ be the first column of $M$ which is not a coordinate vector
(since $M$ is a standard representative, it follows that the entries
in the first $k-1$ column vectors agree with the identity matrix up
to sign).  If we denote the $j^{\mbox{\tiny th}}$ entry of $v_i$ by $v_{ij}$,
then $v_{kk}$ is nonzero.  This follows from the fact
that $M$ is a standard representative and the fact that $M$ is
invertible.

Since the $k^{\mbox{\tiny th}}$ column is not a coordinate vector (and since $M$ is
invertible), there is a $j$, different from $k$, so that $v_{jk}$ is
nonzero.  If there is a $j > k$ such that $v_{jk} \neq 0$, then, since
$M$ is a standard representative, $|v_{jk}| \leq |v_{kk}|$, and $a_k +
\epsilon_j a_j$ is a descending ideal edge for some $\epsilon_j = \pm 1$.  If
$v_{jk}=0$ for all $j > k$, it follows that $v_{kk} = \pm 1$ (since
$M$ is invertible) and there is some $j < k$ so that $v_{jk} \neq 0$.
But then, again, $a_j + \epsilon a_k$ is descending for some $\epsilon
= \pm 1$.
\epf

Here is the proof that descending links are connected:

\begin{proof}[Proof of  Proposition~\ref{proposition:descending
  links are connected}]

We first claim that if $\iota$ is any descending ideal edge, then
there is a subordinate 2-letter ideal edge $\iota'$ which is also
descending; see Section~\ref{section:ideal edges} for definitions.  It
will then follow from Lemma~\ref{lemma:subordinate} and
Corollary~\ref{lemma:description of descending} that there is a path in
$\dlkr$ between the 1-edge blowup realizing $\iota$ to the 1-edge
blowup realizing $\iota'$ (the graph simultaneously realizing $\iota$
and $\iota'$ is the midpoint of the path).

To prove the claim, we need some notation.  First, recall the
notations $\rho$, $a_i$, $v_i$, and $M$ from above.  Also, say
(without loss of generality) that $\iota = a_{i_1} + a_{i_2} + \cdots +
a_{i_m}$, and denote $v_{i_j}$ by $w_j$.  Starting with the matrix
with the $w_i$ as rows, we obtain a matrix $M'$ by deleting all
columns without a nonzero entry.
The $ij^{\mbox{\tiny th}}$ entry
of $M'$ is denoted $w_{ij}$.

We proceed in two cases.  If the first column of $M'$ is not a
coordinate vector, then at least two of the $w_{i1}$ are nonzero, in
particular, $w_{11} \neq 0$.  Without loss of generality, say $w_{11}
> 0$.  Since $\iota$ is descending, there must be a $k$ so that
$w_{k1} < 0$, and since $M$ is a standard representative, we have
$|w_{k1}| \leq |w_{11}|$.  It follows that $a_{i_1} + a_{i_k}$ is
descending, and this completes the proof of the first case.

If the first column of $M'$ is a coordinate vector (i.e. $w_{11} = \pm
1$ and $w_{k1} = 0$ for $k > 1$), then we look at the second column of
$M'$.  Without loss of generality, assume $w_{22} > 0$.  At this point
there are three subcases.  If $w_{k2}=0$ for all $k > 2$, then
$a_{i_1}+a_{i_2}$ is descending, since $\iota$ is descending.  If there is a $k > 2$ so that
$w_{k2} < 0$ then $a_{i_2}+a_{i_k}$ is descending (since $M$ is a
standard representative).  If $w_{k2} \geq 0$ for all $k > 2$ and
$w_{k2} \neq 0$ for at least one $k > 2$, then, since $\iota$ is descending, it follows that $w_{12}
< 0$ and so $a_{i_1}+a_{i_k}$ is
descending for any $k > 2$ with $w_k > 0$.

We now claim that given any two descending 2-letter ideal edges, there
is a path between the corresponding points in $\dlkr$.  This follows, as above, from
Lemma~\ref{lemma:subordinate} and Corollary~\ref{lemma:description of
  descending}, in addition to the fact that opposite 2-letter ideal
edges cannot both be descending (Lemma~\ref{forbidden pairs}).
This completes the proof.
\epf


\subsection{Completely descending link}

We now shift our attention to Proposition~\ref{desc link dim}.  Let
$\rho$ be a rose represented by a marked graph $(\Ga,g)$, and say that
$\Ga$ has edges $a_1, \dots, a_n$ labelled by $v_1, \dots, v_n$.

The main argument for the proof (Section~\ref{section:dimension of
descending links} below) is purely combinatorial, referring only to
isomorphism types of labelled graphs.  As things stand, however, we
cannot describe $\dlkr$ in terms of combinatorial graphs
without metrics.  Indeed, given a marked graph in $\dlkr$,
if we shrink the descending edges to have length less than 1 (while
staying in the frontier by enlarging a nondescending edge), then the
resulting marked graph is not in $\dlkr$
(Corollary~\ref{lemma:description of descending}).

To remedy this problem we perform a deformation retraction of $\dlkr$
onto the \emph{completely descending link}, which we define to be the
subset of $\dlkr$ consisting of marked graphs where each edge not
labelled $\pm v_i$ is descending.  The deformation retraction is
achieved by simply shrinking all edges which correspond to
nondescending ideal edges.  Recall that these edges form a forest (Proposition~\ref{lemma2}(\ref{forest})),
so there is no obstruction.
We denote the completely descending link of $\rho$ by $\cdlkr$.

\begin{lem}
\label{def ret}
For any given rose $\rho$, the completely descending link $\cdlkr$ is
a strong deformation retract of the descending link $\dlkr$.  In particular, the two
are homotopy equivalent.
\end{lem}

We see that $\cdlkr$ has the desired cell structure: a cell
is given by a combinatorial type of labelled graph and the cells are
parameterized by the lengths of the edges in the graph.  To be more
precise, let $(\Ga'g')$ be a marked graph in $\cdlkr$, and for each
$i$, let $k_i$ be the number of edges of $\Ga'$ labelled $\pm v_i$.
For each $i$ we thus get a $(k_i-1)$-simplex by projecting
\[ \set{(t_1, \dots, t_{k_i}) \in [0,1]^{k_i}: t_j=1 \mbox{ for some
  } j} \]
to the simplex $\Delta_i = \set{\sum t_i=1}$.  This projection is a
homeomorphism.  For each edge not labelled $\pm v_i$, we allow its
length to vary arbitrarily within $[0,1]$, as long as one such edge
has length 1.  If $k_0$ is the
number of such edges, then, as above, we get a $(k_0-1)$-simplex
$\Delta_0$.  Thus, the cell corresponding to $(\Ga',g')$ has a cell
structure given by the product:
\[ \Delta_0 \times \cdots \times \Delta_n \]
We now summarize some of the important features of this cell
structure:
\begin{prop}
\label{prop:cdlk}
Consider a cell $C$ of $\cdlkr$ as above.
\begin{enumerate}
\item Passing to faces of $C$ corresponds to collapsing forests in $\Ga'$.
\item $C$ is top-dimensional if and only if all vertices of $\Ga'$
  have valence 3.
\item If $\Ga'$ has $v$ vertices, then $C$ has dimension $v-2$.
\end{enumerate}
\end{prop}


\subsection{Proof of Proposition~\ref{desc link dim}}
\label{section:dimension of descending links}

In this section we show that the completely descending link for any
rose is homotopy equivalent to a complex of dimension $2n-5$ (Proposition~\ref{cdl dim}).  Since
the completely descending link is a deformation retract of the
descending link (Lemma~\ref{def ret}), Proposition~\ref{desc link dim} follows as a corollary.

As usual, let $\rho=(\Ga,g)$ be a rose in $\Y_n$, with edges $a_1,
\dots, a_n$ labelled by $v_1, \dots, v_n$.  If $(\Ga',g')$ is any
marked graph in $\str$, we define the \emph{$v_i$-loop} as the image of $a_i$
under a homotopy inverse of the collapsing map $\Ga' \to \Ga$.

\begin{prop}
\label{cdl dim}
Let $n \geq 3$.  For any rose $\rho$ in $\Y_n$, there is a strong deformation retraction of $\cdlkr$ onto a
complex of dimension $2n-5$.
\end{prop}

\bpf

If any top-dimensional cell of $\cdlkr$ has a free face in
$\cdlkr$, then there is a homotopy equivalence (deformation
retraction) of $\cdlkr$ which collapses away this cell.  We perform
this process inductively until we arrive at a subcomplex $L$ where no top-dimensional cell
has a free face.

We now suppose that $L$ is $(2n-4)$-dimensional,
i.e., it has at least one top-dimensional cell.  Among
these, choose a cell $C$ where the total number of edges
$\ell$ of a $v_1$-loop is minimal.  Call the loop $P$ and choose one
of its edges labelled $\pm v_1$ and call it $e$; see the leftmost
diagram in Figure~\ref{fig:generic}.  Say that $C$ is given
by a marked graph $(\Ga',g')$.

Firstly, note that $\ell$ is not 1, since there are no graphs with
separating edges in $\Y_n$.

\begin{figure}[ht]
\psfrag{v1}{$v_1$}
\includegraphics[scale=.5]{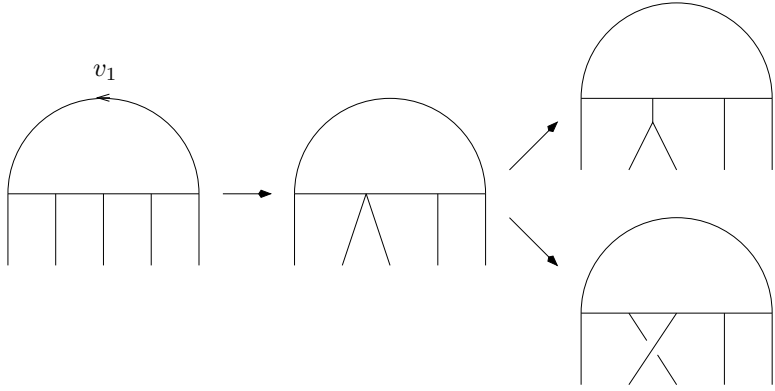}
\caption{The top edge is $e$.  The horizontal path plus $e$ is $P$.}\label{fig:generic}
\end{figure}

If we collapse any edge of $P-e$ (middle of Figure~\ref{fig:generic}),
we move to a codimension 1 face of $C$.  There are two ways to move to
a new top-dimensional cell, since there are two other blowups of the
resulting valence 4 vertex.  One way reduces the length of $P$ (top
right of Figure~\ref{fig:generic}), so by the minimality assumption for
$C$, this is not a cell of $L$.  Since we are assuming $C$ does not
have any free faces, the other top-dimensional cell (bottom right of
Figure~\ref{fig:generic}), call it $C'$, must be in $L$.  The marked
graphs in $C$ and $C'$ have the same labels outside of $P$; the
difference is that the order of the edges leaving $P$ has changed
(Proposition~\ref{lemma1}(\ref{labels dont change}) is applied twice).

Continuing in this way, we see that if we permute the edges leaving
$P$ in any way, we arrive at cells which are necessarily part of $L$.
In particular, the graph obtained by taking the edge which leaves $P$
at one endpoint of $e$ and moving it to the other endpoint of $e$
gives a descending cell $\bar{C}$.

We now argue that $C$ and $\bar{C}$ are opposite in the sense of
Lemma~\ref{forbidden pairs}.  Consider either endpoint of $e$ in
$\Ga'$.  This is a valence 3 vertex, as shown in Figure~\ref{fig:end
of e}.  By Proposition~\ref{lemma1}(\ref{switch condition}) and
Proposition~\ref{lemma2}(\ref{coefficients}), the labels must be as in
the left hand side of the figure.  When we move the edge labelled
$\sum k_iv_i$ to the other end of $e$ (as above), the labels must be
as shown in the right hand side of the figure; the key point is that
the labels and orientations do not change for $e$ and the edge being
moved.  It is then possible for us to determine the label for the
third edge leaving the vertex where these edges meet.  By
Lemma~\ref{forbidden pairs} and Corollary~\ref{lemma:description of
descending}, we have a contradiction.
\epf

\begin{figure}[ht]
\psfrag{1}{$v_1$}
\psfrag{2}{$\displaystyle \sum_{i\geq 2} k_iv_i$}
\psfrag{3}{$v_1+\displaystyle \sum_{i\geq 2} k_iv_i$}
\psfrag{4}{$v_1-\displaystyle \sum_{i\geq 2} k_iv_i$}
\includegraphics[scale=.75]{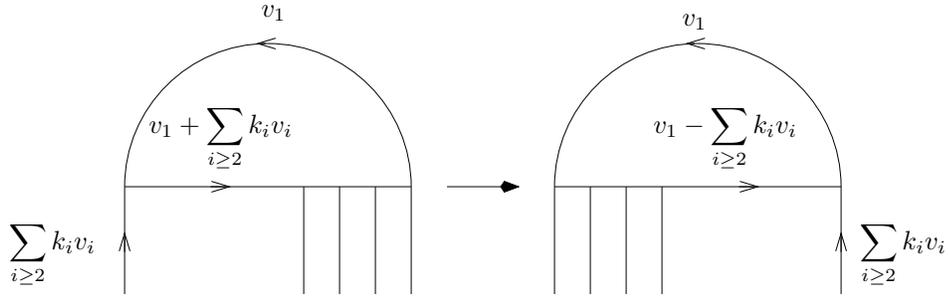}
\caption{Labels at endpoints of $e$.}
\label{fig:end of e}
\end{figure}


\section{The toy model and infinite generation of top homology}
\label{section:infinite generation}
\label{section:toy model}

In this section we prove the second part of the main theorem, that
$H_{2n-4}(\T_n,\Z)$ is not finitely generated when $n \geq 3$.  In order to do this,
we define a subcomplex $\tm_n$ of $\Y_n$, called the ``toy model'', we
find an explicit infinite basis for $H_{2n-4}(\tm_n,\Z)$, and then we
show that the inclusion $\tm_n \to \Y_n$ induces a monomorphism
$H_{2n-4}(\tm_n,\Z) \to H_{2n-4}(\Y_n,\Z)$
(Theorem~\ref{prop:independent tori}).


\subsection{Description of the toy model}
\label{section:description of the toy model}

Let $\rho=(\Ga,g)$ be the rose in $\Y_n$ with the identity marking,
let $x_i$ denote the edges of the standard rose $R_n$, and let $a_i$ denote the corresponding edges of $\Ga$.

Consider the set of points $\tm_n^0=\set{(\Gamma',g')}$ in $\str$ where $g'(x_1) \cup
g'(x_2)$ is a rank 2 rose.  We define the \emph{toy model} to be the
subset $\tm_n$ of $\Y_n$ given by: \[ \tm_n = \bigcup_{p_i,q_i\in \Z} \tm_n^0 \cdot
\left[
\begin{array}{cccccc}
1 & 0 & p_3 & p_4 & \cdots & p_n \\
0 & 1 & q_3 & q_4 & \cdots & q_n \\
0 & 0 & 1& 0 & \cdots & 0 \\
0 & 0 & 0& 1 & \cdots & 0 \\
\vdots & & & & \ddots & \vdots \\
0 & 0 & 0& 0 & \cdots & 1 \\
\end{array}
\right ]
\]

\p{Another point of view}  We now give a different description of
$\tm_n$, which will make it easier to find its homotopy type.

For each marked graph $(\Ga',g')$ of $\tm_n$, the union $g'(x_1) \cup g'(x_2)$
is a rank 2 rose in $\Ga'$, and $\Ga'$ has $n-2$ edges $a_3, \dots,
a_n$ labelled $v_3, \dots, v_n$, where $v_i$ is a coordinate vector
with $+1$ in the $i^{\mbox{\tiny th}}$ spot.
By considering the starting and ending points of $a_3, \dots, a_n$ as
points in $g'(x_1) \cup g'(x_2)$, a path in $\tm_n$ can be thought of
as a path in the configuration space of $n-2$ pairs of points in the
universal abelian cover of $g'(x_1) \cup g'(x_2)$, which is $U = (\R
\times \Z) \cup (\Z \times \R)$.  To make this precise, for each
metric graph $(\Ga',g')$ in $\tm_n$, we rescale the metric so that
$g'(x_1)$ and $g'(x_2)$ both have length 1.  After doing this, the
endpoints of the $a_i$ give a well-defined subset of the metric cover
$U$.

If, in the configuration space, we move the two points
corresponding to the endpoints of some $a_i$ by the same integral
vector, then the corresponding point in $\tm_n$ does not change.

\begin{prop}
\label{proposition:toy model as config space}
The above construction defines a homeomorphism:
\[ (U^2)^{n-2}/(\Z^2)^{n-2} \to \tm_n \]
\end{prop}

At this point, the proof is straightforward and is left to the reader.

A typical graph in $\tm_7$ is shown in Figure~\ref{fig:torus
graph}.  That graph is ``maximally blown up'' in the sense that it has
the greatest number of valence 3 vertices possible in $\tm_7$.

\begin{figure}[ht]
\includegraphics[scale=.5]{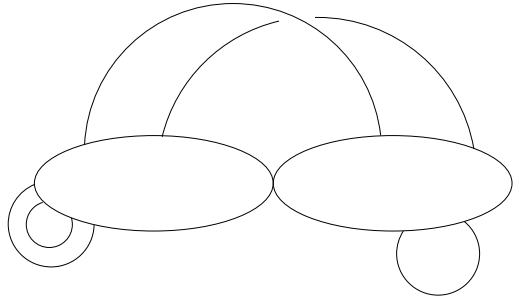}
\caption{A maximally blown up graph in $\tm_7$.}
\label{fig:torus graph}
\end{figure}


\subsection{Homotopy type of the toy model}
\label{section:rank 3 toy model}

We start by focusing our attention on the rank 3 toy model
$\tm_3$.  In general, we have $\tm_n \cong (\tm_3)^{n-2}$, and so we will be able to deduce the finiteness properties of $\tm_n$ from those of $\tm_3$.

Via Proposition~\ref{proposition:toy model as config space}, we can
think of $\tm_3$ as pairs of points in $U$.  However, it will simplify
our analysis if we thicken $U$ to a space $V$, which we now
define.  First, for any integers $p$ and $q$, denote by $D_{p,q}$ the
open disk of radius $r$ around $(p,q) \in \R^2$, for some fixed $r$
close to zero.  Then, define $V = \R^2 - \cup D_{p,q}$.

The straight line retraction of $V$ onto $U$ gives a homotopy
equivalence from $V^2/\Z^2$ to $U^2/\Z^2 \cong \tm_3$.  Thinking of
$V^2/\Z^2$ as pairs of points in $V$, we immediately see the following
features:

\begin{enumerate}
\item The diagonal of $V^2/\Z^2$ is a torus with one boundary
  component.
\item For each $(p,q)$, there is a
  2-torus $Z_{p,q} = \partial D_{0,0} \times \partial D_{p,q}$.
\end{enumerate}

We will now use Morse theory to argue that these features capture
the homotopy type of $\tm_3$.  We consider the Morse function
$d: V^2/\Z^2 \to \R$ which assigns to a point in $V^2/\Z^2$ the
Euclidean distance between the pair of points in $V^2$.

\begin{figure}[ht]
\psfrag{minset}{minset}
\psfrag{index 1}{index 1}
\psfrag{index 2}{index 2}
\includegraphics[scale=.4]{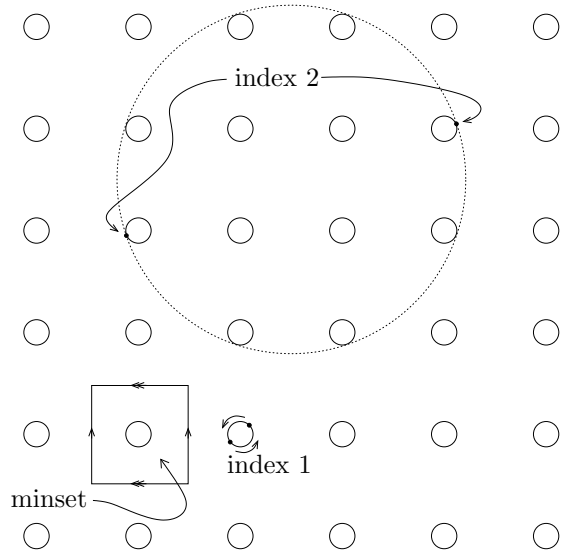}
\caption{Critical points for the toy model in rank 3.\label{critpointspic}}
\end{figure}

We see that $d$ has the following features, depicted in Figure~\ref{critpointspic}:
\blwide

\item The minset is the torus with one boundary component
corresponding to the diagonal of $V^2/\Z^2$.

\item There is a horizontal 1-cell of index 1 critical points
corresponding to pairs of points lying diametrically
opposite from each other on $\partial D_{0,0}$.

\item For every $(p,q)\neq (0,0)$, there is an index 2 critical point,
  corresponding to the two points of tangency of $\partial
  D_{0,0}$ and $\partial D_{p,q}$ with the unique circle tangent to both.

\item At all other points, there is a well-defined gradient flow, and
  so there are no other critical points.
\el

We notice that each index two critical point is the maximum point of
the corresponding $Z_{p,q}$.  Thus, at each critical point of index 2,
a 2-cycle is added.  One can also see that, at the horizontal 1-cell,
there is another torus being added.  Since there are no critical
points of index greater than 2, these classes are nontrivial in $H_2(\tm_3,\Z)$.

In higher rank, we define the torus $Z_{p,q}$ to be the $(n-2)$-fold
product $(\partial D_{0,0} \times \partial D_{p,q})^{n-2}$.
As $\tm_n \cong \tm_3^{n-2}$, we have:

\begin{prop}
\label{proposition:generators for top homology of toy model}
The $Z_{p,q}$ freely generate $H_{2n-4}(\tm_n,\Z)$.
\end{prop}

To formalize the above argument, one can use Morse theory for
manifolds with corners (see \cite{db}).  There is one technicality:
the critical values of $d$ are not isolated; this is easily overcome
by replacing each $D_{p,q}$ with an ellipse (or proving a more
general Morse theory).

We remark that the image of $\pi_1(Z_{p,q}) \cong \Z^{2n-4}$ in
$\pi_1(\Y_n) \cong \T_n$ is a conjugate of the the subgroup $G$
of $\T_n$ described in the introduction.  To see this, one simply
needs to understand the effect of blowups and blowdowns on the
\emph{homotopy} classes of marked graphs; see~\cite{cv}.


\subsection{Independence in homology}  We now set out to prove
that the $Z_{p,q}$ represent independent classes in $\Y_n$
(Theorem~\ref{prop:independent tori}).  In particular, this will
prove the second part of the main theorem.

By understanding the homotopy equivalences
\[ (V^2)^{n-2}/(\Z^2)^{n-2} \to (U^2)^{n-2}/(\Z^2)^{n-2}
\to \tm_n \]
we can give a concrete description of the $Z_{p,q}$ in terms of marked
graphs.  First of all, we have:

\begin{lem}
\label{seg has tori}
Each $Z_{p,q}$ is contained in the union of stars of roses with marking matrices of the form:
\[ \left [
\begin{array}{cccccc}
1 & 0 & p_3 & p_4 & \cdots & p_n \\
0 & 1 & q_3 & q_4 & \cdots & q_n \\
0 & 0 & 1& 0 & \cdots & 0 \\
0 & 0 & 0& 1 & \cdots & 0 \\
\vdots & & & & \ddots & \vdots \\
0 & 0 & 0& 0 & \cdots & 1 \\
\end{array}
\right ]
\]
where each $p_i \in [p-1,p+1]$ and $q_i \in [q-1,q+1]$.
\end{lem}

The main fact we will need about the stars of roses listed in
Lemma~\ref{seg has tori} is that any ideal edge of the form
\[ \pm v_1 \pm v_2 + \sum_{i \geq 3} k_i v_i \]
is \emph{ascending}.
We will also need the following observation about the $Z_{p,q}$:

\begin{lem}
\label{lemma:weakly realized}
For any $Z_{p,q}$ and any rose $\rho$, we can give $Z_{p,q}$ the same cell structure as $\cdlkr$.  In particular, if a marked graph in some $Z_{p,q}$ has an edge of length less than 1 corresponding to an ideal
edge $\iota$, then there is a point in that $Z_{p,q}$ which realizes
$\iota$.
\end{lem}

For the remainder, fix a $Z_{p,q}$, and let $\rho$ be a rose of
greatest norm which intersects $Z_{p,q}$.  By Lemma~\ref{seg has
tori}, if $p$ and $q$ are both nonzero, then $\rho$ is unique; if
one of them is zero, then there are $2^{n-2}$ choices for $\rho$; and if
$p=q=0$, then there are $2^{2n-4}$ choices.
From Lemma~\ref{seg has tori}, we deduce the following key fact:

\begin{lem}
\label{big torus}
Any map from $\set{Z_{p,q}}$ to roses, sending $Z_{p,q}$ to any rose $\rho$
of maximal norm with $Z_{p,q} \cap \str \neq \emptyset$, is injective.
In particular, given any finite subset of $\set{Z_{p,q}}$, the rose of
maximal norm intersecting this set has nonempty intersection with
exactly one torus in this set.
\end{lem}

Let $\dlkr$ denote the descending link for $\rho$ as
defined in Section~\ref{subsection:the induction}.

\begin{lem}
\label{torus sphere}
The intersection $Z_{p,q} \cap \dlkr$ is homeomorphic to $S^{2n-5}$.
\end{lem}

\bpf

We assume $p,q > 0$, with the other cases handled similarly.

Under the homotopy equivalence $(V^2)^{n-2}/(\Z^2)^{n-2} \to
(U^2)^{n-2}/(\Z^2)^{n-2}$, we can identify $Z_{p,q}$ with the
configuration space of $n-2$ pairs of points in $U \subset \R^2$
where the first point $z_i$ in each pair lies on the coordinate square
with vertices at $(0,0)$ and $(1,1)$ and the second point $z_i'$ in
each pair lies on the square with vertices $(p,q)$ and $(p+1,q+1)$.

The rose $\rho$ is realized when each $z_i$ is at $(0,0)$ and each
$z_i'$ is at the point $(p+1,q+1)$.  The points of $Z_{p,q} \cap
\dlkr$ are exactly the set of points where each $z_i$ is within a
distance of $1/2$ from the origin, each $z_i'$ is within $1/2$ of
$(p+1,q+1)$, and at least one $z_i$ or $z_i'$ has distance exactly $1/2$.

In other words, each $z_i$ and $z_i'$ is allowed to move within a
closed interval, and such a configuration is in $Z_{p,q} \cap \dlkr$
if at least one of the points is on the boundary of its interval.
Thus, $Z_{p,q} \cap \dlkr$ is homeomorphic to $\partial I^{2n-4} \cong
S^{2n-5}$.
\epf

The following completes the proof of the main theorem:

\begin{thm}
\label{prop:independent tori}
Let $n\geq 3$.  The $Z_{p,q}$ form an infinite set of independent
classes in $H_{2n-4}(\Y_n,\Z) \cong H_{2n-4}(\T_n,\Z)$.  In other
words, $H_{2n-4}(\tm_n,\Z)$ injects into $H_{2n-4}(\Y_n,\Z)$.
\end{thm}

\bpf

Given a finite subset $A$ of $\set{Z_{p,q}}$, let $Z_{p,q}$ be an element
which intersects a rose $\rho$ of highest norm (marking matrix as in
Lemma~\ref{seg has tori}, the $i^{\mbox{\tiny th}}$ row corresponds to
$v_i$).  We know that there is a strong deformation retraction of
$\dlkr$ onto a complex of dimension $2n-5$ (Lemma~\ref{def ret} plus
Proposition~\ref{cdl dim}).  The goal is to show that the image of the
sphere $Z_{p,q} \cap \dlkr$ is embedded in this complex.  Even better,
we will show that the deformation retractions of Lemma~\ref{def ret}
and Proposition~\ref{cdl dim} do not move the points of $Z_{p,q} \cap
\dlkr$.

To see why this proves this proposition, we consider the long exact
sequence associated to the pair $(\str,\dlkr)$:
\[
\begin{array}{l} \cdots \to H_{2n-4}(\str) \to H_{2n-4}(\str,\dlkr)
 \\  \qquad
\to H_{2n-5}(\dlkr) \to H_{2n-5}(\str) \to \cdots
\end{array}
 \]
By excision, $Z_{p,q}$ corresponds to a class in
$H_{2n-4}(\str,\dlkr)$.  The image in $H_{2n-5}(\dlkr)$ is the class
$Z_{p,q} \cap \dlkr$, which is nontrivial once we show $Z_{p,q} \cap
\dlkr$ is embedded in the
$(2n-5)$-dimensional deformation retract of $\dlkr$ (Lemma~\ref{torus sphere}).  It follows
that  $Z_{p,q}$ is nontrivial in $H_{2n-4}(\str,\dlkr)$ and hence, via
excision, in $H_{2n-4}(\seg',\seg)$ where $\seg$ is the largest initial segment
not containing $\rho$ and $\seg'=\seg\cup\str$.  By Lemma~\ref{big
  torus}, each element of $A$ other than $Z_{p,q}$ is trivial in
$H_{2n-4}(\seg',\seg)$, and so $Z_{p,q}$ is linearly independent from
these, which is what we wanted to show.

Thus, we are reduced to showing that the two deformation retractions do not move the sphere $Z_{p,q} \cap \dlkr$.  We handle each in turn.

For the deformation retraction of $\dlkr$ onto $\cdlkr$, we need to
show that $Z_{p,q} \cap \dlkr$ is already contained in $\cdlkr$.
Suppose that there were a point of $Z_{p,q} \cap \dlkr$ which were not
contained in $\cdlkr$.  By
Lemma~\ref{lemma:weakly realized}, there is a point which realizes an
\emph{ascending} ideal edge (Lemma~\ref{lemma:neighbors}), and this implies that $Z_{p,q}$ intersects some rose
of higher norm, contradicting the choice of $\rho$.

We now focus on the deformation retraction of Proposition~\ref{cdl
dim}.  A marked graph representing a maximal cell of $\cdlkr$ must
have disjoint $v_1$- and $v_2$-loops.  Indeed, there are no valence 4 vertices, so the overlap would have to contain an ascending edge by
Proposition~\ref{lemma1}(\ref{switch condition}), the statement after
Lemma~\ref{seg has tori}, and Lemma~\ref{descending criterion}.  The
codimension 1 cells which get collapsed are obtained from these
maximal cells by collapsing an edge of the $v_1$-loop.  Thus, the
corresponding graphs still have disjoint $v_1$- and $v_2$-loops.  On
the other hand, in any graph of $Z_{p,q}$, the $v_1$-loop and the
$v_2$-loop intersect in exactly 1 point.  Thus, no points of $Z_{p,q}$
are moved during this retraction, so we are done.
\epf


\appendix

\section*{Appendix: Proof of Proposition~\ref{magnus}}
\label{section:magnus}

This appendix contains Magnus's proof of Proposition~\ref{magnus}.  In
this section, we freely use the notation of
Section~\ref{section:finite generation}.

Let $K$ be the subgroup of $\out(F_n)$ generated by the $K_{ik}$ and
$K_{ikl}$ for distinct $i$, $k$, and $l$.  We now prove
Proposition~\ref{magnus}, that $K$ is normal in $\out(F_n)$.

\bpf[Proof of Proposition~\ref{magnus}]

We choose the following generating set for $\out(F_n)$:
\[ \begin{array}{rcl}
\delta_{12} & = & [x_1 x_2, x_2, \dots, x_n] \\
\Omega_1 & = & [x_1^{-1}, x_2, \dots, x_n] \\
\Pi_{i-1} & = & [x_1, \dots, x_{i-2}, x_i, x_{i-1}, x_{i+1},
    \dots, x_n]
\end{array}\]

It suffices to show that the conjugates of the $K_{ik}$ and
$K_{ikl}$ by the chosen generators of $\out(F_n)$ (and their
inverses) are elements of $K$.

We have the following simplifications:

\begin{enumerate}

\item Operations on disjoint sets of elements commute.

\item Since $\Omega_1$ and $\Pi_{i-1}$ have order 2, we don't need to
conjugate by their inverses.

\item We don't need to conjugate by $\delta_{12}^{-1}$ since \[ (\Pi_1
\Omega_1 \Pi_1) \delta_{12} (\Pi_1 \Omega_1 \Pi_1)^{-1} =
\delta_{12}^{-1} \]

\item Since $K_{ikl} = K_{ilk}^{-1}$, we may assume $k < l$.

\item  Any outer automorphism $\psi$ of the form
\[ [x_1, \dots, g' x_i g, \dots, x_n] \]
where $gg'$ is an element of the commutator subgroup of the subgroup $H$
of $F_n$ generated by $\set{x_k:k \neq i}$ is an element of $K$.

\bigskip

To see that $\psi \in K$, first note that, by postcomposing with a product of $K_{i\star}^{\pm 1}$,
we may assume that $g'=1$.  Now, we know that the commutator subgroup of $H$ is normally generated
by the $[x_k,x_l]$, where $k$ and $l$ are both different from $i$.
Therefore, it suffices to handle the case of
\[ g = h [x_k,x_l] h^{-1} = [hx_kh^{-1},hx_lh^{-1}] \]
where $h = x_{i_1} \cdots x_{i_p}$ is an arbitrary element of $H$.
It is elementary to check that
\[ \psi = P^{-1} K_{ikl} P \]
where
\[ P = \prod_{j \neq i} K_{ji_p} \cdots \prod_{j \neq i} K_{ji_2}
\prod_{j \neq i} K_{ji_1} \]
\end{enumerate}

\bigskip

Given these simplifications, it is straightforward to check (case by case) that the
conjugates by $\Pi_{i-1}$, $\Omega_1$, and $\delta_{12}$ of each
$K_{ik}$ and $K_{ikl}$ are elements of  $K$.  There is one exception; we give Magnus's
computation for this difficult case here:
\[ \delta_{12} K_{2l1} \delta_{12}^{-1} =
K_{l2}K_{l1}^{-1}K_{1l}K_{l1}K_{2l1}K_{12l}K_{l2}^{-1}K_{2l}^{-1} \]
\epf

\bibliographystyle{plain}
\bibliography{tn}

\begin{thebibliography}{10}

\bibitem{rk}
Problems in low-dimensional topology.
\newblock In Rob Kirby, editor, {\em Geometric topology (Athens, GA, 1993)},
  volume~2 of {\em AMS/IP Stud. Adv. Math.}, pages 35--473. Amer. Math. Soc.,
  Providence, RI, 1997.

\bibitem{bt}
Gilbert Baumslag and Tekla Taylor.
\newblock The centre of groups with one defining relator.
\newblock {\em Math. Ann.}, 175:315--319, 1968.

\bibitem{db}
Dietrich Braess.
\newblock Morse-{T}heorie f\"ur berandete {M}annigfaltigkeiten.
\newblock {\em Math. Ann.}, 208:133--148, 1974.

\bibitem{cv}
Marc Culler and Karen Vogtmann.
\newblock Moduli of graphs and automorphisms of free groups.
\newblock {\em Invent. Math.}, 84(1):91--119, 1986.

\bibitem{ah}
Allen Hatcher.
\newblock {\em Algebraic topology}.
\newblock Cambridge University Press, Cambridge, 2002.

\bibitem{dj1}
Dennis Johnson.
\newblock The structure of the {T}orelli group. {I}. {A} finite set of
  generators for {${\mathcal I}$}.
\newblock {\em Ann. of Math. (2)}, 118(3):423--442, 1983.

\bibitem{km}
Sava Krsti{\'c} and James McCool.
\newblock The non-finite presentability of {${\rm IA}(F\sb 3)$} and {${\rm
  GL}\sb 2({\bf Z}[t,t\sp {-1}])$}.
\newblock {\em Invent. Math.}, 129(3):595--606, 1997.

\bibitem{wm}
W.~Magnus.
\newblock {\"Uber n-dimensional Gittertransformationen}.
\newblock {\em Acta Math.}, 64:353--367, 1934.

\bibitem{mcm}
Darryl McCullough and Andy Miller.
\newblock The genus {$2$} {T}orelli group is not finitely generated.
\newblock {\em Topology Appl.}, 22(1):43--49, 1986.

\bibitem{gm}
Geoffrey Mess.
\newblock The {T}orelli groups for genus {$2$} and {$3$} surfaces.
\newblock {\em Topology}, 31(4):775--790, 1992.

\bibitem{jn}
J.~Nielsen.
\newblock {Die isomorphismengruppe der freien Gruppen}.
\newblock {\em Math. Ann.}, 91:169--209, 1924.

\bibitem{sv1}
John Smillie and Karen Vogtmann.
\newblock Automorphisms of graphs, {$p$}-subgroups of {${\rm Out}(F\sb n)$} and
  the {E}uler characteristic of {${\rm Out}(F\sb n)$}.
\newblock {\em J. Pure Appl. Algebra}, 49(1-2):187--200, 1987.

\bibitem{sv2}
John Smillie and Karen Vogtmann.
\newblock A generating function for the {E}uler characteristic of {${\rm
  Out}(F\sb n)$}.
\newblock In {\em Proceedings of the Northwestern conference on cohomology of
  groups (Evanston, Ill., 1985)}, volume~44, pages 329--348, 1987.

\bibitem{kv}
Karen Vogtmann.
\newblock Automorphisms of free groups and outer space.
\newblock In {\em Proceedings of the Conference on Geometric and Combinatorial
  Group Theory, Part I (Haifa, 2000)}, volume~94, pages 1--31, 2002.

\end{thebibliography}

\end{document}